\documentclass{icota}

\usepackage[dvips]{graphicx}

\usepackage{amssymb}
\usepackage{amsfonts}
\usepackage{amsmath}

\newtheorem{theo}{Theorem}[section]
\newtheorem{coro}{Corollary}[section]
\newtheorem{lemma}{Lemma}[section]
\newtheorem{defin}{Definition}[section]
\newtheorem{prop}{Proposition}[section]
\newtheorem{examp}{Example}[section]

\begin{document}

\title{\Large \bf G-coupling functions:
The Infinite Dimensional case}

\author{
{\bf \underline{D.M. Morales Silva}}\\
{\it School of ITMS}, {\it University of Ballarat}, {\it VIC
3353}, {\it Australia}\\ dmoralessilva@students.ballarat.edu.au }
%\and
%{\bf Author B}\\
%{\it Affiliation B}, {\it University B}, {\it Address B}, {\it Country B}\\
%nameB@authorB }

\maketitle

\maketitle
\thispagestyle{empty}

\begin{center}
{\bf Abstract}
\end{center}
In this work we present a class of functions, motivated by gap
functions, which we call G-coupling functions. We will show that
these functions can generate a duality scheme for minimization
problems by means of the general conjugation theory. Thanks to
this scheme, a Lagrange-type function is introduced as well.

\section{Introduction}

For solving non-convex optimization problems, a tool that is
becoming more important is \emph{generalized conjugation}. In
\cite{MOR} the G-coupling functions are introduced in finite
dimensional spaces. Here we extend this definition to the infinite
dimensional case. These coupling functions will allow us to see
duality schemes in a different way. The usual theory found in the
literature (\cite{RUB}, \cite{SOSA}, and references therein) are
related to a fixed coupling function, but here we consider (for a
specified function $f$) a family of coupling
functions.\\

These coupling functions are motivated by gap functions. It is
interesting to point out, that many of these (gap) functions have
similar properties. However, in some cases they are functions of
one vector and it is important, since they are linked to specified
optimization problems, that those functions have zeros.\\

On the other hand, G-coupling functions will be defined as
functions in two variables and they might not have zeros. Even
more, given a specified proper function $f$, it is shown that a
certain sub-family of this family of coupling functions satisfies
many interesting properties.\\

In Section 2, we describe how many gap functions have similar
properties, which are useful for the definition of G-coupling
functions.\\

In Section 3, it is found the definition of G-coupling function
with properties related to generalized conjugation using this
family of functions and a fixed proper function $f$.\\

In Section 4, it can be seen how these ideas generate
Lagrange-type functions (see \cite{RUB-YANG}).

\section{Motivation}

In several works already published, there can be found definitions
of GAP functions for particular problems. Now we present 2
concrete examples.\\

In \cite{BUR}, the Variational Inequality Problem is studied:
$$(VIP)\text{ Find }x_0\in C,\text{ such that, }\exists y^*\in
T(x_0)\text{ with }\langle y^*,x-x_0\rangle \geq 0\ \forall x\in
C,$$ where $T$ is a maximal monotone correspondence which is
defined as follows: given a point to set map, $T$, it will be said
that it is a maximal monotone correspondence if it satisfies that
$\langle u-v,x-y \rangle \geq 0$ for every $u\in T(x),\ v\in T(y)$
with $x,y\ \in C$ and if there exists $v$, such that $\langle
u-v,x-y \rangle \geq 0$, for all $x,y\ \in C$ and for all $u\in
T(x)$, then $v\in T(y)$. The corresponding GAP function is then
defined as follows:
$$h_{T,C}(x):=\sup_{(v,y)\in G_C(T)}\langle v,x-y\rangle,$$ where
$G_C(T)=\{(v,y):\ v\in T(y),\ y\in C\}$ and $C$ is a non-empty
closed convex set. This function happens to be non-negative and
convex, and it is equal to zero only in solutions
of $(VIP)$.\\

In \cite{SOSA}, the Equilibrium Problem is studied:
$$(EP)\text{ Find }x\in K,\text{ such that }f(x,y)\geq 0,\ \forall
y\in K,$$ where $K\subset I\!\!R^n$ is a non-empty closed convex
set and $f:K\times K\rightarrow I\!\!R$ is a function that
satisfies:
\begin{enumerate}
\item[i)] $f(x,x)=0$, for all $x\in K$.
\item[ii)] $f(x,\cdot):K\rightarrow I\!\!R$ is convex and l.s.c.
\item[iii)] $f(\cdot,y):K\rightarrow I\!\!R$ is u.s.c.
\end{enumerate}
The GAP function is defined as:
$$g_f(y):=\left\{
\begin{array}{cc}
   \displaystyle\sup_{x\in K}f(x,y) & if\ y\in K  \\
   +\infty & \text{in other case}. \\
\end{array}
\right.$$ In this case, the function $g_f$ is non-negative, convex
and l.s.c. and if it vanishes at $x_0$, then $x_0$ is a solution
of $(EP)$.\\

%In \cite{YANG}, the Extended Pre-Variational Inequality Problem is
%studied:
%$$(EPVIP)\text{ Find }x_0\in I\!\!R^n,\text{ such that }\langle
%F(x_0),\eta(x,x_0)\rangle \geq f(x_0)-f(x),\ \forall x\in
%I\!\!R^n,$$ where $F:I\!\!R^n \rightarrow I\!\!R^n$, $\eta:
%I\!\!R^n\times I\!\!R^n \rightarrow I\!\!R^n$ and $f:I\!\!R^n
%\rightarrow I\!\!R\cup\{+\infty\}$. In this work, the GAP function
%is $$\min_{y\in I\!\!R^n}[\langle F(x),\eta(y,x)\rangle
%-f(x)+f(y)],$$ which is non-positive and it only reaches the value
%zero in solutions of $(EPVIP)$.\\

In these examples, gap functions are used to transform a special
Equilibrium Problem (for example, the VIP is a particular case of
an EP) into a minimization problem.\\

Now our attention is focused in using coupling functions that
could be related, at least in some general aspect, to GAP
functions. Therefore these functions must link both primal and
dual variables. Since these coupling functions must be related to
a sense of ``gap", we consider these functions as non-negative
and with 2 arguments.\\

Let us remember that for the minimization problem, the convex
conjugation theory allows us to generate a dual problem and there
is implicit another concept of gap function (see \cite{AVRIEL},
\cite{J.P.SOSA.OC} and \cite{ROCK}): consider
$$\alpha =\inf[f(x):x\in I\!\!R^n]. \qquad (P)$$ Define a function
$\varphi:I\!\!R^n\times I\!\!R^p\rightarrow \overline{I\!\!R}$,
where $\overline{I\!\!R}=I\!\!R\cup\{-\infty,+\infty\}$,
satisfying
$$\varphi(x,0)=f(x),\ \forall x\in I\!\!R^n.$$ Then $\varphi$ will
be called a perturbation function and the function
$h:I\!\!R^p\rightarrow \overline{I\!\!R}$ defined by
$$h(u)=\inf_{x\in I\!\!R^n}\varphi(x,u)$$ will be called the
marginal function. Observe that $$\alpha=h(0)=\inf_{x\in
I\!\!R^n}\varphi(x,0)=\inf_{x\in I\!\!R^n}f(x).$$ Considering now
$h^{**}$, the convex bi-conjugate (see \cite{ROCK}) of $h$ one
has:
$$h^{**}(0)\leq h(0)=\alpha$$ where $$h^{**}(0)=\sup[\langle
u^*,0\rangle -h^*(u^*): u^*\in I\!\!R^p].$$ Then, making
$-\beta=h^{**}(0)$, one has $$\beta=\inf_{u^*\in
I\!\!R^p}h^*(u^*).\qquad (Q)$$ $(Q)$ is called the dual problem of
$(P)$ and in general we have $-\beta\leq\alpha$. It is said that
there is no duality gap whenever $h^{**}(0)=h(0)$. It is easy to
prove that $h^*(u^*)=\varphi^*(0,u^*)$, and if we define the
function $k:I\!\!R^n\rightarrow \overline{I\!\!R}$ by
$\displaystyle k(x^*):=\inf_{u^*\in I\!\!R^p}\varphi^*(x^*,u^*)$,
then $\beta=k(0)$.\\

This analysis is summarized in the following scheme:
$$\begin{array}{rclcrcl}
  \alpha & = & \inf f(x) \qquad (P) & & \beta & = & \inf h^*(u^*) \qquad (Q) \\
  \varphi(x,0) & = & f(x),\ \forall x\in I\!\!R^n & & \varphi^*(0,u^*) & = & h^*(u^*),\ \forall u^*\in I\!\!R^p \\
  h(u) & = &\displaystyle \inf_x \varphi(x,u) & & k(x^*) & = &\displaystyle \inf_{u^*}\varphi^*(x^*,u^*) \\
  \alpha & = & h(0) & & \beta & = & k(0) \\
\end{array} $$
$$-\beta \leq \alpha.$$
If $h$ is proper and convex, a necessary and sufficient condition
for ensuring that there will be no duality gap ($-\beta =\alpha$)
is that $h$ be l.s.c. at 0 (in general $\varphi$ l.s.c. does
not imply that $h$ would be l.s.c.).\\
\\
Further more, if $h$ is convex, l.s.c. and $0\in ri(dom(h))$, then
$\alpha=-\beta$ and the dual problem has at least one optimal
solution, and if $\overline{u^*}$ is an optimal solution of $(Q)$
and $\varphi=\varphi^{**}$, then
$$\overline{x} \text{ is an optimal solution of }(P)
\Longleftrightarrow f(\overline{x})+h^*(\overline{u^*})=0.$$
Consider now the function $g:I\!\!R^n \times I\!\!R^p \rightarrow
\overline{I\!\!R}$ defined by: $$g(x,u^*)=f(x)+h^*(u^*).$$ This
function vanishes at $(x_0,u^*_0)$ if and only if $x_0$ solves the
primal problem and $u^*_0$ solves the dual one. In addition, this
function is non-negative and if the first variable is kept fixed,
the function is convex and l.s.c. It is clear now, which
properties are satisfied for many gap functions.

\section{G-coupling Functions}

As stated before, G-coupling functions are first introduced in
\cite{MOR} for finite dimensional spaces. We are going to extend
this notion for arbitrary Banach spaces.\\

Henceforth, we consider two arbitrary Banach spaces $X$ and $Y$.
\begin{defin}
A non-negative function $g:A\times B\rightarrow I\!\!R$, with
$A\times B\subset X\times Y$ will be called a G-coupling function
if
\begin{enumerate}
\item[(D1)] $\displaystyle \inf_{x\in A,\ y\in B}g(x,y)=0.$

%\item[iii)] $g(x,\cdot):C \rightarrow I\!\!R$ is a
%convex and l.s.c. function for each $x$ in $I\!\!R^n$.
\end{enumerate}
\end{defin}

Define
\begin{equation}\mathcal{F}^{A,B}:=\{g:A\times B\rightarrow
I\!\!R:\ g \text{ is a G-coupling function}\}.\end{equation}

Not every G-coupling function has zeros:\\ \\ {\bf Example:}
Define on $X\times Y$ $$g(x,y)=\exp(\|x\|-\|y\|).$$ Then
$g\in\mathcal{F}^{X,Y}$ is continuous and it does not have
any zeros.\\

Let us turn our attention now to how the family of functions
$\mathcal{F}^{A,B}$ will allow us to establish duality schemes in
(at least for now) the minimization problem. It is important to
point out that in the following we consider an unusual type of
duality, $f:A\rightarrow I\!\!R\cup\{+\infty\}$ is kept fixed and
$g\in \mathcal{F}^{A,B}$, for a given $B\subset Y$, is variable.\\

Consider a proper function $f:A\rightarrow I\!\!R\cup\{+\infty\}$.
For a given $B\subset Y$ take $g\in \mathcal{F}^{A,B}$. Define
$f^g:B\rightarrow I\!\!R\cup\{+\infty\}$ and $f^{gg}:A\rightarrow
I\!\!R\cup\{+\infty\}$ as follows (for example see \cite{RUB} and
references therein): \begin{equation} f^g(y):=\sup_{x\in
A}\{g(x,y)-f(x)\}\ \forall y\in B,\end{equation}
\begin{equation}f^{gg}(x):=\sup_{y\in B}\{g(x,y)-f^g(y)\}\ \forall
x\in A.\end{equation} In some cases, it would be better to
consider a $g\in \mathcal{F}^{A,B}$ which
satisfies:\emph{\begin{enumerate}
\item[(D2)] $B$ is convex and $g(x,\cdot):B\rightarrow I\!\!R$ is a
convex and l.s.c. function for each $x$ in $A$.
\end{enumerate}} With this, we have the following:

\begin{lemma} Let $f:A\rightarrow I\!\!R\cup \{+\infty\}$ be a
proper function and given $B\subset Y$ take $g\in
\mathcal{F}^{A,B}$. Then
\begin{equation}f^g(y)+f(x)\geq g(x,y)\geq 0,\ \forall (x,y)\in
A\times B,\end{equation} which implies
\begin{equation} f(x)\geq -f^g(y),\ \forall (x,y)\in A\times B.
\end{equation}
Moreover if $g$ satisfies $(D2)$, then $f^g$ is a convex l.s.c
function.
\end{lemma}

Unless it is mentioned, not every $g\in \mathcal{F}^{A,B}$
satisfies $(D2)$.\\

It would be interesting to know which condition either a
G-coupling function $g$ or the function $f$ must satisfy in order
that the function $f^g$ be proper, because with this one would
have a non-trivial function related to $f$. The following lemma
ensures the existence of such a function $g\in \mathcal{F}^{A,B}$
for any $B\subset Y$, taking as a starting point a natural
condition on $f$ which must be imposed if $f$ is the objective
function of a minimization problem.

\begin{lemma} Let $f$ be as before. Then $f$ is bounded from below
if and only if, for every non-empty $B\subset Y$, there exists
$g\in \mathcal{F}^{A,B}$ such that $f^g$ is proper.
\end{lemma}
{\bf Proof:}\begin{enumerate}\item[$\bullet$] Suppose that $\inf
f>-\infty$, then for a non-empty $B_0\subset Y$ fixed, consider
$g\in \mathcal{F}^{A,B_0}$ as follows:
$$\displaystyle g(x,y)= \|y\|,\ \forall (x,y)\in A\times B.$$ Thus
$$f^{g}(y)=\|y\|-\inf f\ \forall y\in B_0,$$
which is clearly a proper function and since $B_0\subset Y$ was
fixed arbitrarily, the result is satisfied for every $B\subset Y$.
\item[$\bullet$] Take a non-empty $B_0\subset Y$ and $g\in
\mathcal{F}^{A,B_0}$ such that $f^g$ is proper. Let us suppose
that $\inf f=-\infty$, from \cite{RUB} we can see that this
implies that $\inf f^{gg}=-\infty$. Then:
$$-\infty=\inf f^{gg}=\inf_{x\in A}\left(\sup_{y\in B_0}
[g(x,y)-f^g(y)]\right)\geq$$ $$\sup_{y\in B_0}\left(\inf_{x\in A}
[g(x,y)-f^g(y)]\right)\geq\sup_{y\in B_0}(-f^g(y))=-\inf_{y\in
B_0}f^g(y),$$ which means $\displaystyle -\infty\geq -\inf_{y\in
B_0}(f^g(y))$. Then $\displaystyle \inf_{y\in B_0}f^g(y)=+\infty$,
which implies that $f^g$ is not proper and we have a
contradiction. Therefore we must have that $\inf f>-\infty$.
\end{enumerate}
Notice that this proof also states, in particular, that there
exists $g\in \mathcal{F}^{A,B}$ for every non-empty $B\subset Y$
which satisfies $(D2)$ and $f^g$ is proper.\\

%Henceforth, consider only functions $f$ such that $\inf f>-\infty$
%and for some fixed $B\subset Y$, $g\in \mathcal{F}^{A,B}$ will be
%such that $f^g$ is proper.\\

Given non-empty sets $A\subset X$ and $B\subset Y$, let
\begin{equation}\mathcal{F}^{A}:=\{f:A\rightarrow I\!\!R\cup\{+\infty\},\
f\text{ is proper, } \inf f>-\infty\}\end{equation} and
$\gamma_{g,f}:A\times B\rightarrow I\!\!R\cup\{+\infty\}$ defined
by: \begin{equation}\gamma_{g,f}(x,y):=f(x)+f^g(y)\end{equation}
with $g\in \mathcal{F}^{A,B}$ and $f\in \mathcal{F}^A$. Take
$f\in\mathcal{F}^A$ and define
\begin{equation}\mathcal{F}_f^{A,B}:=\{g\in \mathcal{F}^{A,B}/ f^g \text{ is
proper and }\inf \gamma_{g,f}=0\}.\end{equation} {\bf Remark:}
Observe that $\gamma_{g,f}$ might not be in $\mathcal{F}^{A,B}$,
since $\gamma_{g,f}$ can take the value $+\infty$ for somewhere in
$A\times B$.
\begin{lemma} $\mathcal{F}_f^{A,B}$ is non-empty for all non-empty
$B\subset Y$.
\end{lemma}
{\bf Proof:} Given a non-empty $B\subset Y$, define $g\in
\mathcal{F}^{A,B}$ by:
$$g(x,y)=\|y\|.$$ It is easy to check that $g$ belongs to
$\mathcal{F}_f^{A,B}$ (this example also proves that there are
functions in $\mathcal{F}_f^{A,B}$ which satisfy $(D2)$).\\

Now consider \begin{equation}(P)\ \min_x f(x)\end{equation} with
$f\in \mathcal{F}^A$. Taking $g\in \mathcal{F}_f^{A,B}$, define
the dual problem related to $g$:
\begin{equation}(D_g)\ \min_{y\in B}f^g(y).\end{equation} Since
$$\inf_{(x,y)\in A\times B}\gamma_{g,f}(x,y)=\inf_{x\in A}
f(x)+\inf_{y\in B}f^g(y)=0,$$ then
\begin{equation}\inf_{x\in A} f(x)=-\inf_{y\in B}f^g(x^*)=
\sup_{y\in B}[-f^g(y)].\end{equation} This means that there is no
duality gap between the primal problem $(P)$ and its dual $(D_g)$
for every $g\in \mathcal{F}_f^{A,B}$.
\begin{theo} Let $g\in\mathcal{F}^{A,B}_f$. Then $\overline{y}$ is
a solution of $(D_g)$ and $\overline{x}$ is a solution of $(P)$ if
and only if $\gamma_{g,f}(\overline{x},\overline{y})=0$.
\end{theo}
{\bf Proof:} $\overline{x}$ and $\overline{y}$ are solutions of
$(P)$ and $(D_g)$ respectively if and only if
$$f(\overline{x})= \inf f=-\inf f^g=-f^g(\overline{y})
\Longleftrightarrow f(\overline{x})+f^g(\overline{y})=\gamma_{g,f}
(\overline{x},\overline{y})=0._{\Box}$$\\ \\ {\bf Remark:} The
previous result suggest us that the function $\gamma_{g,f}$ can be
seen as the GAP function of problem $(P)$ and its dual $(D_g)$.\\

The next theorem states that given non-empty sets $A\subset X$ and
$B\subset Y$, the correspondence defined by
$$\begin{array}{cccc}
  \mathbf{F}: & \mathcal{F}^A & \rightrightarrows &
  \mathcal{F}^{A,B}\\
   & f & \mapsto & \mathbf{F}(f)=\mathcal{F}_f^{A,B},
\end{array}$$ is a closed correspondence (see \cite{ZANG}).

\begin{theo}
Take $f\in \mathcal{F}^A$ ($A\subset X$ is non-empty) and a
non-empty $B\subset Y$. If there exist $f_k:dom(f)\rightarrow
I\!\!R$, $g_k:A\times B\rightarrow I\!\!R$, sequences of functions
($k\in I\!\!N$), such that:
\begin{enumerate}
\item[i)] $f_k$ converges uniformly to $f$ in $dom(f)$.% $$\forall
%\varepsilon>0,\exists k_0\in I\!\!N\text{ such that
%}|f_k(x)-f(x)|\leq \varepsilon, \forall k\geq k_0\text{ and }x\in
%dom(f).$$
\item[ii)] $g_k\in \mathcal{F}_{f_k}^{A,B}$ satisfies $(D2)$ for every
$k\in I\!\!N$.
\item[iii)] $g_k$ converges uniformly to a function $g$ in
$A\times B$.%: $$\forall \varepsilon>0,\exists k_1\in I\!\!N\text{
%such that }|g_k(x,y)-g(x,y)|\leq \varepsilon, \forall k\geq
%k_1\text{ and }(x,y)\in A\times B.$$
\end{enumerate} Then $g\in \mathcal{F}_f^{A,B}$ and it
satisfies $(D2)$.
\end{theo}
{\bf Proof:} Let us prove first that $g\in \mathcal{F}^{A,B}$.
Since $g_k$ converges uniformly to $g$, given $\varepsilon>0$,
there exists $N\in I\!\!N$ such that if $k\geq N$ then
$$|g_k(x,y)-g(x,y)|< \varepsilon,\qquad \forall (x,y)\in
A\times B.$$ $$\text{Hence }
g_k(x,y)-\varepsilon<g(x,y)<g_k(x,y)+\varepsilon,\ \forall
(x,y)\in A\times B.$$ Taking $\displaystyle \inf_{x,y}$ (remember
that $\inf g_k =0$ for all $k\in I\!\!N$):
$$-\varepsilon<\inf_{x,y}g(x,y)< \varepsilon.$$ Then $|\inf
g|<\varepsilon$. And since $\varepsilon>0$ is arbitrary, one has
that $\inf g=0$. This proves that $g\in\mathcal{F}^{A,B}$.\\
\\ Now we prove that $g$ satisfies $(D2)$. We need to prove that
$g(x,\cdot):B\rightarrow I\!\!R$ is convex and l.s.c. for all
$x\in A$. Let $x_0\in A$ be fixed arbitrarily.
\begin{enumerate}
\item[$\bullet$] $g(x_0,\cdot)$ is convex: since for all $k\in
I\!\!N$, $g_k(x_0,\cdot)$ is convex, one has that given
$y_1,y_2\in B$ and $t\in [0,1]$:
$$g_k(x_0,ty_1 + (1-t)y_2)\leq
tg_k(x_0,y_1)+(1-t)g_k(x_0,y_2).$$ Making $k\rightarrow +\infty$:
$$g(x_0,ty_1 + (1-t)y_2)\leq tg(x_0,y_1)+(1-t)g(x_0,y_2),$$
which proves that $g(x_0,\cdot)$ is convex.
\item[$\bullet$] $g(x_0,\cdot)$ is l.s.c.: fix $y_0\in B$ and take
$\lambda <g(x_0,y_0)$. There exists $N\in I\!\!N$ such that
$$|g_N(x,y)-g(x,y)|<\varepsilon,\ \forall (x,y)\in A\times B,$$
where $\displaystyle \varepsilon= \frac{g(x_0,y_0)-\lambda}{2}$.
$$\text{Hence } \lambda< \lambda +
\varepsilon=g(x_0,y_0)-\varepsilon<g_N(x_0,y_0).$$ Since
$g_N(x_0,\cdot)$ is l.s.c., then there exists $V(y_0)\subset B$, a
neighborhood of $y_0$, such that if $y\in V(y_0)$ then
$$\lambda+\varepsilon<g_N(x_0,y).$$ Reducing $g(x_0,y)$:
$$\lambda+\varepsilon-g(x_0,y)<g_N(x_0,y)- g(x_0,y)<\varepsilon.$$
Therefore, if $y\in V(y_0)$, then $\lambda<g(x_0,y)$. Thus
$g(x_0,\cdot)$ is l.s.c. in $y_0\in B$, and since $y_0$ was fixed
arbitrarily then $g(x_0,\cdot)$ is a l.s.c. function.
\end{enumerate}
We have proved that for a fixed $x_0\in A$, $g(x_0,\cdot)$ is a
convex l.s.c. function, and since $x_0$ was fixed arbitrarily we
have proved in fact that $g\in\mathcal{F}^{A,B}$ satisfies $(D2)$.\\
\\ It remains to prove that $g\in \mathcal{F}_f^{A,B}$. For
doing this, let us show that $(f_k^{g_k})_{k\in I\!\!N}$
converges uniformly to $f^g$ (in $B$).\\
\\ Let $\varepsilon>0$ and $N\in I\!\!N$ be such that if $k\geq N$
then $$|g_k(x,y)-g(x,y)|<\frac{\varepsilon}{4},\ \forall (x,y)\in
A\times B$$ and $$|f_k(x)-f(x)|<\frac{\varepsilon}{4},\ \forall
x\in dom(f).$$ Fix $k\geq N$ and take $y\in B$ arbitrarily, then
$$f_k^{g_k}(y)-\frac{\varepsilon}{2} < g_k(x',y)-f_k(x'), \text{
for some }x'\in dom(f).$$ Hence
$$f_k^{g_k}(y)-\varepsilon<g_k(x',y)-f_k(x')
-\frac{\varepsilon}{2}<g(x',y)-f(x')\leq f^g(y),$$ and so
\begin{eqnarray}\label{stabeq1}
f_k^{g_k}(y)-\varepsilon<f^g(y).
\end{eqnarray} This proves that
$f_k^{g_k}(y)-f^g(y)<\varepsilon.$ On the other hand:
$$f^g(y)-\frac{\varepsilon}{2}<g(x'',y)-f(x''),\ \text{for some
}x''\in dom(f),$$ whence $$f^g(y)-\varepsilon
<g(x'',y)-f(x'')-\frac{\varepsilon}{2}<g_k(x'',y)-f_k(x'')\leq
f_k^{g_k}(y),$$ and so $$f^g(y)-\varepsilon< f_k^{g_k}(y).$$ This
shows that
\begin{eqnarray}\label{stabeq2}
-\varepsilon<f_k^{g_k}(y)-f^g(y).
\end{eqnarray}
Since $y\in B$ was fixed arbitrarily, thanks to (\ref{stabeq1})
and (\ref{stabeq2}) we have that
$$-\varepsilon<f_k^{g_k}(y)-f^g(y)<\varepsilon,\text{ for every
} y\in B.$$ This proves that $(f_k^{g_k})_{k\in I\!\!N}$ converges
uniformly to $f^g$ (in $B$), and it is immediate to see that $f^g$
is proper and
$$0\leq f(x)+f^g(y)\leq f_k(x)+f_k^{g_k}(y)+\varepsilon,\ \forall
(x,y)\in dom(f)\times B,$$ where $\varepsilon>0$ is arbitrary and
$k$ is large enough. Taking $\displaystyle\inf_{(x,y)\in A\times
B}$ one has: $$0\leq \inf_{(x,y)\in A\times B}(f(x)+f^g(y))\leq
\varepsilon.$$ Therefore $\displaystyle\inf_{(x,y)\in A\times
B}(f(x)+f^g(y))=0$
and $g\in \mathcal{F}_f^{n,m}._\Box$\\
\\ This theorem proves a more difficult situation, the case when
$g_k\in \mathcal{F}_{f_k}^{A,B}$ satisfy $(D2)$ for all $k\in
I\!\!N$. For the general case, just omit the two $\bullet$ items
and change $B$ for a non-empty set.\\

At this point a natural question arises, for given
$f\in\mathcal{F}^A$ and $g\in \mathcal{F}_f^{A,B}$, would be there
any kind of relation between the optimal points and the optimal
values of $f$ and $f^{gg}$? The next lemma answers this.
\begin{lemma}
For a fixed non-empty $B\subset Y$ and every $g\in
\mathcal{F}^{A,B}_f$, the following are satisfied:
\begin{enumerate}
\item[i)] $\inf f=\inf f^{gg},$
\item[ii)] if $x_0$ is a global minimum of $f$, then $x_0$ is a
global minimum of $f^{gg}$.
\end{enumerate}
\end{lemma}
{\bf Proof:} Remember that $f^{gg}$ is defined by:
$$f^{gg}(x)= \sup_{y\in B}\{g(x,y)-f^g(y)\}.$$
\begin{enumerate}
\item[i)] $\inf f^{gg}\leq \inf f$ is always true. On the other
hand $$f^g(y)+f^{gg}(x)\geq g(x,y)\geq 0,\ \forall x\in A,\ y\in
B,$$ which implies that $$\inf_{x\in A} f^{gg}(x)\geq -\inf_{y\in
B}f^g(y).$$ But, since $g\in \mathcal{F}_f^{A,B}$ one has that
$$\inf f=-\inf_{y\in B}f^g(y),$$ which means
$$\inf f\leq \inf f^{gg}\leq \inf f.$$ Therefore $\inf f=\inf
f^{gg}$.
\item[ii)] $f^{gg}(x_0)\leq f(x_0)=\inf f=\inf f^{gg}\leq
f^{gg}(x_0),$ then $f^{gg}(x_0)=\inf f^{gg}$.
\end{enumerate}

%\section{Generalized Lagrangians}

\section{Lagrangians induced by $\mathcal{F}_f^{A,B}$}

Take $f\in \mathcal{F}^A$, a non-empty $B\subset Y$, $g\in
\mathcal{F}_f^{A,B}$ and consider $$(P):\ \inf_{x\in A}f(x).$$
Recall that
$$(D_g):\ \min_{y\in B}f^g(y)$$ is the dual problem of $(P)$
related to $g$. Define $L_1:I\!\!R^n \times C\rightarrow
I\!\!R\cup\{+\infty\}$, as follows:
\begin{equation}L_1(x,y):=f(x)-g(x,y).\end{equation} This function
has some interesting properties:

\begin{theo}
\begin{equation}\sup_{y\in B}\inf_{x\in A} L_1(x,y)=\inf_{x\in
A} \sup_{y\in B}L_1(x,y).\end{equation}
\end{theo}
{\bf Proof:} The inequality $\displaystyle \sup_{y\in B}\inf_{x\in
A} L_1(x,y)\leq\inf_{x\in A} \sup_{y\in B}L_1(x,y)$ is always
true. For the opposite:
$$L_1(x,y)=f(x)-g(x,y)\leq f(x),\ \forall (x,y)\in
A\times B,$$ then
$$\sup_{y\in B}L_1(x,y)\leq f(x),\ \forall x\in A.$$ It
follows that $$\inf_{x\in A}\sup_{y\in B}L_1(x,y)\leq \inf_{x\in
A} f(x).$$ But, since $g\in \mathcal{F}_f^{A,B}$, we have that
$$\inf_{x\in A} f(x)=-\inf_{y\in B}f^g(y)=-\left(\inf_{y\in
B}\left\{ \sup_{x\in A} [g(x,x^*)-f(x)] \right\} \right)$$
$$\Longrightarrow \inf_{x\in A} f(x)= \sup_{y\in B}\inf_{x\in A}
L_1(x,y),$$ which means, $$\inf_{x\in A}\sup_{y\in B}L_1(x,y)\leq
\sup_{y\in B}\inf_{x\in A} L_1(x,y).$$ Finally,
$$\sup_{y\in B}\inf_{x\in A} L_1(x,y)=\inf_{x\in A} \sup_{y\in
B}L_1(x,y)._\Box$$

We are interested now in which properties are satisfied for every
saddle-point of $L_1$. Remember that $(x_0,y_0)\in A\times B$ is a
saddle point of $L_1$ if and only if $$ L_1(x_0,y)\leq
L_1(x_0,y_0)\leq L_1(x,y_0),\ \forall (x,y)\in A\times B.$$

\begin{prop}\label{propsaddleGENERAL}
Let $L_1$ be as before, if there exists $(x_0,y_0)\in A\times B$
saddle point of $L_1$, then:
\begin{enumerate}
\item[i)] $x_0 \in dom(f)$.
\item[ii)] $y_0$ is an optimal solution of $(D_g)$.
\item[iii)] $f^{gg}(x_0)=f(x_0)$.
\end{enumerate}
\end{prop}
{\bf Proof:}\begin{enumerate}
\item[i)] This is immediate thanks to the definition of saddle point.
\item[ii)] From the previous theorem and the definition of saddle
point, we have that
$$L_1(x_0,y_0)= \sup_{y\in B}\inf_{x\in A}
L_1(x,y)=\inf_{x\in A} \sup_{y\in B}L_1(x,y).$$ But
$$\sup_{y\in B}\inf_{x\in A} L_1(x,y)=-\inf_{y\in
B}f^g(y),$$ moreover
$$L_1(x_0,y_0)=\inf_{x\in A} L_1(x,y_0)=-f^g(y_0).$$
Thus, $$f^g(y_0)=\inf_{y\in B}f^g(y).$$
\item[iii)] $\displaystyle f^{gg}(x_0)=\sup_{y\in B}[g(x_0,y)-f^g(y)]=
\sup_{y\in B}\left[g(x_0,y)-\sup_{z\in A} [g(z,y)-f(z)]\right].$
Which means, $$f^{gg}(x_0)=\sup_{y\in B}\inf_{z\in
A}[g(x_0,y)-g(z,y)+f(z)]= \sup_{y\in B}\inf_{z\in
A}[g(x_0,y)+L_1(z,y)].$$ This implies $$f^{gg}(x_0)\geq \inf_{z\in
A}[g(x_0,y_0)+L_1(z,y_0)]=g(x_0,y_0)+\inf_{z\in A} L_1(z,y_0),$$
but since $(x_0,y_0)$ is a saddle point of $L_1$, then
$\displaystyle\inf_{z\in A} L_1(z,y_0)=L_1(x_0,y_0)$. With this,
we have that
$$f^{gg}(x_0)\geq g(x_0,y_0)+L_1(x_0,y_0)=f(x_0),$$ which means
$f^{gg}(x_0)\geq f(x_0)$. $f^{gg}(x_0)\leq f(x_0)$ is always true
(see \cite{RUB} and references therein).$_\Box$
\end{enumerate}

\begin{prop}
If $x_0$ is a solution of $(P)$ and $x_0^*$ is a solution of
$(D_g)$, then $(x_0,x_0^*)$ is a saddle point of $L_1$.
\end{prop}
{\bf Proof:} This is immediate from $$0\leq g(x_0,x_0^*)\leq
f(x_0)+f^g(x_0^*)=0\Longrightarrow g(x_0,x_0^*)=0._\Box$$

In Proposition \ref{propsaddleGENERAL} we would like to improve
the fact that, in general, for every saddle point $(x_0,y_0)\in
A\times B$ of $L_1$ we have that $f^{gg}(x_0)=f(x_0)$. For doing
this, we impose an additional condition over $g$.

\begin{prop} Let $g\in \mathcal{F}^{A,B}_f$ be such that
$\displaystyle \inf_{y\in B} g(x,y)=0$ for every $x\in A$. The
following are equivalent:
\begin{enumerate}
\item[i)] $(x_0,y_0)$ is a saddle-point of $L$.
\item[ii)] $x_0$ is a solution of $(P)$ and $y_0$ is a solution
of $(D_g)$.
\end{enumerate}
\end{prop}
{\bf Proof:} The implication ii) $\Rightarrow$ i) is true thanks
to the previous Proposition.\\ \\ Consider now $(x_0,y_0)$ a
saddle-point of $L_1$, then $$L_1(x_0,y)\leq L_1(x_0,y_0),\
\forall y\in B,$$ which is equivalent to $$f(x_0)-g(x_0,y)\leq
f(x_0)- g(x_0,y_0),\ \forall y\in B$$ $$\Updownarrow$$
$$g(x_0,y_0)\leq g(x_0,y),\ \forall y\in B.$$ Finally
$$g(x_0,y_0)=\inf_{y\in B}g(x_0,y)=0.$$ On the other hand
$$L_1(x_0,y_0)\leq L_1(x,y_0),\ \forall x\in A.$$ This implies
that $$f(x_0)\leq f(x)-g(x,y_0),\ \forall x\in A$$ (remember that
$g(x_0,y_0)=0$). Taking $\displaystyle \inf_{x\in A}$ we have
$$f(x_0)\leq -f^g(y_0).$$ And thus $f(x_0)=-f^g(y_0)$, which
means that $x_0$ is a solution of $(P)$ and $y_0$ is a solution of
$(D_g)$.\\
\\ {\bf Remark:} To prove that there exists a $g\in
\mathcal{F}^{A,B}_f$ such that $\displaystyle \inf_{y\in
B}g(x,y)=0$ for every $x\in A$ just consider the trivial function
$g\equiv 0$.

\section*{Examples}

For these examples, consider $X=I\!\!R^n$, $h:X\rightarrow
I\!\!R^m$, $$A:=\{x\in X:h(x)\leq 0\}$$ and $f:A\rightarrow
I\!\!R$.

\begin{enumerate}
\item[\bf 1.]{\bf Classical Lagrangian}
%\subsection{Classical Lagrangian}

Let $Y=I\!\!R^m$ and $h$ be such that $h_i:I\!\!R^n\rightarrow
I\!\!R$ is convex and l.s.c. for all $i=1,\ldots,m$. Consider
$$(CP):\ \min_{x\in A}f(x),$$ where $f$ is convex and l.s.c.\\

Remember that (see \cite{AVRIEL} and \cite{J.P.SOSA.OC}) the
following is the well known dual problem:

$$(D_L):\ \min_{\lambda^*\geq 0}\sup_{x\in A}\{\langle
\lambda^*,-h(x)\rangle -f(x)\},$$ $h(x)=(h_1(x),\ldots, h_m(x))$.
Moreover, $x_0$ is a solution of $(CP)$ and $\lambda^*_0$ is a
solution of $(D_L)$ if and only if $(x_0,\lambda^*_0)$ is a saddle
point of the Lagrangian function $L$, given by
$$L(x,\lambda^*):=f(x)+\langle \lambda^*,h(x)\rangle,\ x\in A,\
\lambda^*\in I\!\!R^m_+.$$

Taking $B:=I\!\!R^m_+$, define $g:A\times B\rightarrow I\!\!R$ as
follows: \begin{equation}g(x,\lambda^*):=\langle
\lambda^*,-h(x)\rangle.\end{equation} It is not difficult to show
that $\mathcal{F}^{A,B}_f$ and, even more,
$$f^g(\lambda^*)=\sup_{x\in A}\{\langle \lambda^*,-h(x)\rangle
-f(x)\},\ \lambda^*\in B.$$ Therefore, using G-coupling functions,
we have recovered the classical lagrangian duality.

\item[\bf 2.] {\bf Non-linear lagrangian function}
%\subsection{Lagrange-type functions}

In \cite{RUB-YANG} we find the following well studied case of a
non-linear \emph{lagrange-type} function:
$$L(x,\omega)=f(x)+\max\{\langle\omega_0,h(x)\rangle,\ldots,
\langle\omega_p,h(x)\rangle\},$$ where $x\in I\!\!R^n$ and
$\omega\in (I\!\!R^m_+)^{1+p}$ ($p\in I\!\!N$).\\

If we consider $Y=(I\!\!R^m)^{1+p}$ and $B=(I\!\!R^m_+)^{1+p}$,
define $g:A\times B\rightarrow I\!\!R$ as follows:

\begin{equation}g(x,\omega):= \min(\langle
-h(x),\omega_0\rangle,\ldots,\langle
-h(x),\omega_p\rangle),\ x\in A,\ \omega\in B,\end{equation} we
will have that $g\in\mathcal{F}^{A,B}_f$ and the lagrangian
function induced is the same Lagrange-type function given by
\cite{RUB-YANG}.
\end{enumerate} {\bf Acknowledgements:} I wish to thank to Dr.
Wilfredo Sosa who was my supervisor during my Master thesis, most
of these ideas are from that work. Special thanks to Dr. Regina
Burachik, Dr. David Yost and Prof. Alex Rubinov for all their
support and comments.

%The following examples show us that the definition of G-coupling
%function can be used not only for inducing dual problems.

%%{\bf Remark:} It is not that obvious that the techniques applied
%%in \cite{RUB-YANG} can still be applied in this setting of
%%G-coupling functions. This considerations will be studied in a
%%future work.
%\begin{enumerate}
%\item[\bf 3.] {\bf Min-type functions}
%%\subsection{Min-type functions}

%Consider now $X=Y=I\!\!R^n$, $A=B=I\!\!R^n_+$ and the following
%function: $$\langle l,x\rangle^-:=\left\{\begin{array}{cl}
%\displaystyle\min_{i\in I_+(l)}l_i x_i, & I_+(l)\neq \emptyset \\
%0, & I_+(l)=\emptyset
%\end{array} \right.$$ with $(x,l)\in A\times B$ and
%$I_+(l):=\{i\in\{1,\ldots,n\}:l_i\neq 0\}$.\\

%In \cite{RUB} this function is studied in an Abstract Convexity
%frame. Now we can affirm that $g(x,l):=\langle l,x\rangle^-$
%belongs to $\mathcal{F}^{A,B}$.

%\item[\bf 4.] {\bf Monotonic Analysis over Cones}
%%\subsection{Increasing convex along rays functions}

%In \cite{MONAN1} and \cite{MONAN2}, the following function is
%considered $l:C\times C\setminus\{0\}\rightarrow I\!\!R$

%$$l(x,y):=\max\{\lambda>0:\lambda y\leq_K x\},$$ where $C\subset
%K\subset X$, $C$ is a cone, $K$ is a closed pointed convex cone
%and $X$ is a locally convex Hausdorff topological vector space.\\

%Here if $Y=X$, $A=C$ and $B=C\setminus\{0\}$, we have that $l\in
%\mathcal{F}^{A,B}$.

%\end{enumerate}

\end{document}